# Notes for a study of the didactic transposition of mathematical proof

Nicolas Balacheff*

It is nowadays common to consider that proof must be part of the learning of mathematics from Kindergarten to University[1]. As it is easy to observe, looking back to the history of mathematical curricula, this has not always been the case either because following an old pedagogical tradition of rote learning proof was reduced to the formalism of a text and deprived from its meaning or, despite its acknowledged presence anywhere in mathematics, proof did not get the status of something to learn for what it is. On the long way from its absence as such in the past to its contemporary presence as a content to be taught at all grades, proof has had to go through a process of didactical transposition to satisfy a number of different constraints either of an epistemic, didactical, logical or mathematical nature. I will follow a chronological order to outline the main features of this process with the objective to better understand the didactical problem that our current research is facing.

## A note on didactic transposition

### The concept

As for any other content, teaching and learning mathematics do require, as far as possible, its complete and precise specification as a knowledge to be taught, be it of a conceptual, methodological or technical nature (e.g. integers, multiplying two binomial – FOIL, factorization, integration of a function). It must first be identified, then uttered and finally chosen. This is the social responsibility of several organizations which include those who put it in use as well as decision makers. A complex process takes charge of its formalization as a content to be taught which outcomes are curricula, textbooks and other texts published and disseminated by the organizations in charge. It is recognized since the origin of modern education that the knowledge in use and the knowledge to be taught share similarities while showing significant differences. These phenomena are conceptualized and modelled by the concept of didactic transposition coined by Yves Chevallard (1985). It refers "to the transformations an object or a body of knowledge undergoes from the moment it is produced, put into use, selected and designed to be taught, until it is actually taught in a given educational institution." (Chevallard & Bosch, 2014, p. 170).

Numerous organizations contribute to the decision to turn a piece of knowledge into a content to be taught and participate in shaping its transposition: professional mathematicians, various users from engineers to dealers, associations of teachers, bodies at the different layers of the educational institution and educational decision makers. The forces which interact are of a professional, political and social nature. Yves Chevallard approaches the related phenomena through "the study of the conditions enabling and the constraints hindering the production, development and diffusion of knowledge and, more generally, of any kind of human activity in social institutions" (ibid. p. 173). The didactic transposition is a human and social enterprise which study is now part of the *Anthropological Theory of Didactic* (Chevallard, 1998; Chevallard et al., 2015). In this theoretical framework, taking into consideration "transpositive phenomena means moving away from the classroom and being provided with notions and elements to describe the bodies of knowledge and practices involved in the different institutions at different moments of time." (Chevallard & Bosch, 2014, p. 172).

---

* Team MeTAH – Laboratoire d'informatique de Grenoble – Université Grenoble Alpes, CNRS, Grenoble INP – Grenoble, France

[1] c.f. (Stylianides & Harel, 2018)

## About the method

The process of didactic transposition is multi-stage, it involves several forms of a piece of knowledge, whose relations are shaped by the interactions between its different contributors. Four forms play a pivotal role: the *scholarly knowledge*, the *knowledge to be taught*, the *taught knowledge* and the *learnt knowledge* (Chevallard & Bosch, 2014, p. 170). Since the objective of these notes is to outline the dynamic of the didactic transposition of mathematical proof from an historical perspective, I will align my quest on the first two forms, the scholarly knowledge and the knowledge to be taught. This will set a limit to this exploration since the distance between the knowledge to be taught and the knowledge taught may be important, as we know that it is the case between the knowledge taught and the learned knowledge. Moreover, there are little historical data on what happens in classroom either from the teaching or the learning perspective. But this limit is not too severe given my purpose. I expect that the chosen approach for these notes informs us sufficiently precisely about the nature and the complexity of the phenomena underlying the transposition process.

These notes are based on published research on the history of mathematics education, especially the history of the teaching of geometry, readings of original treatises and textbooks, institutional comments and curricula. The study of these resources is faced to *the fragility of knowledge* which paradoxically despite a "constant form" does not keep constant meaning (Kang & Kilpatrick, 1992, p. 3). This classic remark for linguistics and communication sciences has profound consequences for the epistemological analysis of processes dedicated to knowledge communication and dissemination since they assume hypothetical students, teachers and classrooms whose perception and models are rarely documented (ibid. p. 5).

Hence, readers must keep in mind that the analysis and comments shared in these notes have a conjectural stance, though they paradoxically pretend a sufficiently robust contribution to research on proof as a content to be taught. This robustness comes from the quality of the data and resources, and their significance in relation to the phenomena identified.

**Caveat**: I use the French quotation marks « … » when inserting a quote translated by myself from the French, and classical quotation marks "…" when the quote is originally in English.

## Geometry, the theoretical ground of proof

### Euclid's elements

Mathematical proof appeared very late as an explicit content to be taught, when considering its early formalization by the Greeks, in the 3$^{rd}$ century BC. It comes nowadays with a vision of writing mathematically with absolute rigour *deductive reasoning* based on *explicit foundations*, definitions and postulates. Indeed, this is an idealization of what underpins *Euclid elements of geometry* which shaped the construction of the concept of proof and stimulated the development of *mathematics as a discipline*. This perennial reference has contributed to the constitution of standards of communication between mathematicians. Even if it has been contested and transformed in the course of the history of mathematics, it constitutes a landmark of its epistemology.

Euclid's Elements are considered as one of the treasures of the ancient Greek legacy, though they have been left out for a long time, until the 12$^{th}$ century. Some geometry was probably taught for its practical inputs, but there are not so many evidences about this period when education was not organized in a systematic way; the quadrivium had neither teachers nor students, and other domains than mathematics had priority (grammar, rhetoric) (Høyrup, 2014, sect. 2 & 3). The first milestone of importance from the perspective of mathematics education, following Jens Høyrup is, somewhere in the first half of the 12$^{th}$ century, the translation of Euclid's Elements from the Arabic, "presumably" by Adelard of Bath (fl. 1116-1142).

It is still too commonly claimed that the journey from Euclid BC to Euclid AD in the West only took a diversion through Arab countries. On the contrary, the Arabic translation of the Elements in the 8th century was the subject of discussion among mathematicians of the time about the text and the usefulness of studying it. In particular, Sonja Brentjes (2014) refers to Al-Sizjī disagreeing "vehemently" with those who downplayed the study of Euclid, and insisting on "the need to study and work hard to become a good geometer" (ibid. p. 92). Al-Sizjī gives reason for the study of the Elements among which:

> "[to follow] the methods of them (these theorems and preliminaries) in a profound and successful way, so that you rely not only on the theorems and preliminaries and constructions and arrangements which we mentioned. But you must combine with that (your own) cleverness and guesswork and tricks. The pivotal factor in this art is the application of tricks, and not only (your own) intelligence, but also the thought of the experienced (mathematicians), the skilled, those who use tricks" (Brentjes, 2014, p. 92).

Then, the 10th century Islamic position was that "the beginner must learn the first theorems of Euclid's Elements" (ibid.).[2]

Later versions of Euclid Elements derived from the first Adelard translation, possibly due to his students (Høyrup, 2014, p. 114), were marked by "didactical concerns": "at this point it is clear that the matter presented in the work had become the primary aim, while further utility for astronomy (and, still further, for astrology) had retreated into the background" (ibid.).

To a certain extent, because it was not compulsory, the Elements have been taught since then but it was a challenge given the scarce educational material and the style of the lectures. This was the *pre-Gutemberg era* (ibid. 119). The print of the Elements in 1482 in Venice radically changed access to the text. Controversies over the Elements and their translations developed until the Latin text of Commandin (1572) which seems to have reached some form of consensus (Loget, 2004).

The interest for Geometry[3] grown during the 16th century together with critical readings of Euclid's Elements for its logic, principles and order of propositions (E. Barbin & Menghini, 2014, p. 474). The perspective was theoretical, questioning the text, hence the Euclidean standard pattern of proof.

### The criticisms of Euclid

Let us open this section with a quotation of the foreword of the *Elemens d'Euclide* of R. P. Dechalles de la *Compagnie de Jésus* (1660), first translated in French by Ozanam, member of the Royal Academy of Sciences, in 1677:

> « Having noticed for a long time, that most of those who learn the elements of Euclid, are very often disgusted with them, for not knowing what is the use of such apparently insignificant propositions, and nevertheless so useful; I thought it would be very appropriate, not only to make them as easy as possible; but also, to add some small uses, after each proposition, which would show their usefulness. This is what obliged me to change some demonstrations, which I judged to be too awkward, and beyond the ordinary reach of

---

[2] Although difficult because of the methodological problems of working with a large historical corpus and avoiding cultural bias in the analysis (Brentjes, 2019), the study of Islamic mathematics education in this period is most definitely needed.

[3] I will write Geometry with a capital letter to refer to pure geometry, a body of theoretical knowledge, as opposed to practical geometry, a body of knowledge meant to be used to solve technical, professional or everyday problems. The latter is not well defined, but one may understand that this kind of knowledge claims some direct usefulness. The boarder between theoretical and practical is not in the form of the text, for example Descriptive geometry has a utilitarian raison d'être and is based on a theoretical construction. Both were taught, and are still taught, but with radically different educational objectives and with different pedagogical approaches. Geometry was the natural terrain of the debates on the teaching of mathematical proof. This distinction may be debatable under the light of specific cases, and I thank the first readers of this note for their comments, but it is useful for guiding the analysis.

> beginners, and to substitute a few more intelligible ones. » (Dechalles, 1660, sect. Avant propos - my translation)

This foreword to the Elements sets out the main issues that would be disputed in the centuries to come: the conflict between theoretical and practical views of geometry[4], and the necessity and complexity of proofs, with Euclid standing out as the key reference for learning mathematics.

The theoretical character of Geometry, which gives priority to the logical organisation of the text, was denounced by Descartes. In his essay on the *Rules for the direction of the mind* (ca 1628)[5], he criticizes the ancient mathematicians (in fact, the Elements) for the priority given to *convincing* to the detriment of *understanding*, and thus their inability to allow the reason and meaning of proofs, and therefore of theorems, to be grasped; that is to say the understanding of Geometry itself. This criticism points to a tension between convincing and explaining that we know now is inherent in argumentative discourse. Descartes conceptualised this tension by distinguishing *analysis* and *synthesis*, the former "is the best and truest method of instruction", while le latter "is very suitable to deploy in geometry"[6]. As a matter of fact, he sets the problem of the teaching of geometry in a way which is still relevant: how to manage the equilibrium between convincing and explaining. These criticisms of the Euclid's Elements had a translation in the writing of the *Géométrie* of Arnaud (1667) who proposed a new organisation following the principles of analysis. Sylvestre-François Lacroix, more than a century later, wrote about Arnaud's Geometry: « *This work, is, I believe, the first one in which the order of the proposals of Geometry agrees with that of abstractions, considering first the properties of lines, then that of surfaces, and then that of bodies*. » (Lacroix, 1799, p. xix).

In the first half of the 18th century, Étienne Bonnot de Condillac in *Essay on the origins of human knowledge* (1746) in which he acknowledge the model of rigour of the Geometry, but criticizing its lack of, let say, simplicity. His fundamental postulate claims the superiority of metaphysics for the formation of « a luminous, precise and extensive mind, and which, consequently, must prepare for the study of all the other sciences » (ibid. p. 13). While he takes geometry as a model for the construction of his essay, he criticises it for its failure to find « an order simple enough and easy enough to arrive at the obvious » (ibid.). His essay includes considerations on the communication of knowledge which had a long-lasting influence: "*Finally, after having developed the progress of the operations of the soul and those of language, I try to indicate the means by which one can avoid error, and to show the order one must follow, either to make discoveries, or* to instruct others *of those one has made*." (ibid. p.16).

In the same period the mathematician Alexis Claude Clairaut published his *Elements de Géométrie* (1741) with a preface which takes a clear didactical position rooted in a pedagogical observation: « [...] *it's common for beginners to get tired & put off, before they have any distinct idea of what we wanted to teach them.* » (ibid. p. ij). Clairaut searched for an approach « [...] *bringing together the two advantages of interesting & enlightening the beginners* ». (ibid. p. iij). Looking into the history of Geometry he chosen the problems of « *measuring lands* » (fr: *arpentage*) to give meaning and to avoid proofs for the obvious because those interested in geometry « *enjoyed exercising their minds a little; & on the contrary, they were put off, when they were overwhelmed with demonstrations, which were almost useless* »[7] (ibid. p. x) but he warns the reader: « *this is no ordinary treatise on land-*

---

[4] This tension is still unsolved. The two conflicting questions "how to cater for the elite" versus "how to cater for the wider group of students for whom mathematics should be grounded in real world problem solving and daily life applications", as expressed by Gert Schubring (2015) in a critical analysis of the "Mathematics for all" movement, have not yet receive proper responses.

[5] Original reading *Règles pour la direction de l'esprit* (ca 1628/1953, p. 37-118)

[6] Original reading *Méditations, objections et réponses* (1641/1953, p. 387 sqq.). English translations of the quotes from (Cunning, 2015)

[7] Ceux qui s'intéressent à la géométrie "se plaisaient à exercer un peu leur esprit ; & au contraire, ils se rebutaient, lorsqu'on les accablait de démonstrations, pour ainsi dire inutiles"

*surveying »* (ibid. p. xij). However, it is not a textbook either. The Elements of Clairaut should be read as a manifesto bringing a contribution to the ongoing discussions on the Euclid Elements (Glaeser, 1983).

At the end of the 18th century, Euclid's Elements are the reference which criticisms nourish the search for a text of Geometry for both its theoretical establishment and its communication. This search focussed on the organisation of the text and the need for proof for some propositions, the Euclidean standard remaining stable.

## The 19th century, an epistemological dispute

The French organisation of higher education at the end of the 18th century widened the gap between the theoretical and the utilitarian approach of geometry. Two different courses of study developed, corresponding to two different systems. On the one hand, a teaching of geometry oriented towards applications, essentially in engineers and military schools, on the other hand, a teaching of Geometry centred on geometry for itself, taught in the Normal schools[8] and Faculties. In this context, the criticisms of Euclid lead not only to the writing of new Elements, but also to the writing of treatises with the objective to satisfy practical needs. A significant example is the contribution of the mathematician Gaspard Monge[9] to the project for a Public instruction:

> « *There is an order of knowledge of an indispensable necessity for the stone masons, stonemasons, carpenters, joiners, carpenters, stone locksmiths, contractors of all kinds, painters, engineers of the bridges and roads, engineering officers [...]. The order of knowledge in question here is based on a particular geometry of the three dimensions which does not exist of a well-made treatise; on a purely descriptive, but rigorous geometry, and the purpose of which is to represent by drawings that have only two dimensions of the objects that have three.* » (Monge 1793 quoted by Eveline Barbin (2021, p. 104)).

This « particular geometry » is the Descriptive geometry, which is not a new geometry but a reliable, robust and efficient geometrical instrument designed on and with Geometry in order to manipulate graphical representations of geometrical objects modelling objects from the physical world. Rigor is evoked and required, but proof is not the central preoccupation. This orientation is very strong, it could be seen as a seed of what is now known as Applied Mathematics. It developed in France with the *Grandes Écoles* system for higher education. It had weaker links with general education and the secondary school system than Universities, which were more theoretically oriented. I will continue this exploration focusing on the latter where mathematical proof as a teaching object emerged at the end of the century.

The 19th century saw the nurturing of educational systems at national scales for primary education and secondary education as well, though the population having access to the latter was rather limited. First, secondary education was expensive, second, it was essentially oriented towards University education, third, it was reserved to boys. Two textbooks played a distinctive role because of their impact in France and abroad: the *Éléments de géométrie* (1794) of Adrien Marie Legendre and the *Éléments de géométrie à l'usage de* l'*École Centrale des quatre nations* (1799) written by Sylvestre-François Lacroix (E. Barbin & Menghini, 2014, sect. 4.1 & 4.2). They influenced teacher education and the teaching of geometry in secondary schools, with multiple editions along the century. The writing of each of these two books was driven by different objectives.

---

[8] "Normal school" is the equivalent of the "College of education" or "Teacher training school" of the contemporary educational systems.

[9] Gaspard Monge is the renowned creator of Descriptive geometry, this course given at the École Normale and École Polytechnique was published in 1799 after notes on his 1795 lectures (his assistant Sylvestre-François Lacroix had published in 1795 an "Essais de géométrie sur les plans et les surfaces courbes – ou "Élémens de géométrie descriptive").

Adrien Marie Legendre acknowledged writing his Elements following the « method » of Euclid and Archimedes. He wrote in his preface: « *in trying to equal or even surpass my models of accuracy, I also wanted to spare the reader as much trouble as I could, and I made my efforts to give the demonstrations all the clarity and brevity that the subject matter entails* » (Legendre, 1794, p. vj-viij). He did not hesitate to use Algebra when relevant because « *it would be childish to always use a laborious method when it can be replaced by a much simpler and safer one.* »[10] (ibid.). The Legendre *Elements,* are of a theoretical nature. This is illustrated by his refusal to use of the fifth of the Euclid postulates, hence ever and ever searching a new proof of the theorem on the sum of the angles of a triangle.

Sylvestre-François Lacroix inherits from Condillac, and from his apprenticeship as an assistant of Gaspard Monge. The long preface to his book—a « *preliminary discourse* » —advocates a priority of analysis over synthesis, understanding over convincing and argues the educational value of Geometry  « *[which] is perhaps, of all the parts of mathematics, the one that one should learn first; it seems to me very likely to interest children, as long as it is presented to them mainly in relation to its applications, either on paper or in the field*. » (1799/1804 p.xxix)[11] However, while its structure is different, the body of the book is written according to the traditional Euclidean style.

Both books fix the geometric terms at the beginning of the text, following the Euclidean tradition. Then, they add second order vocabulary defining the terms axiom, theorem, corollary and problem. With respect to the Greek tradition, this is an innovation[12].

The *Elements of Legendre* (1794) includes in the second order terms a definition of theorem and demonstration in the same sentence: « *Theorem is a truth that becomes evident by means of a reasoning called demonstration.* » (ibid. p. 4). The text is organised as a sequence of books, themselves made of a sequence of numbered Propositions. For each Proposition, a subtitle indicates its type, Theorem or Problem, then comes the text in italics of the theorem or of the problem immediately followed by the text of its proof or its solution in roman character.  An important part of Legendre's Elements are notes on some of the demonstrations.

The *Elements of Lacroix* (1798), defines the mathematical terms and second order terms. He does not define the word demonstration, but indicate that *a theorem is a statement which must be demonstrated*.  The structure of a theorem is made explicit indicating that it has two parts, a hypothesis and a conclusion and warning the reader that their role cannot in general be exchanged (ibid. p. lxxxviij). Then the presentation of the text deviates from the Euclidean model. It is organised in two parts, themselves cut into sections with subsections entitled following its theme. In a section, comes a sequence of numbered subsections which subdivisions are labelled first indicating the nature of the subsection (theorem, problem) then the nature of the related discourse (resp. demonstration, solution). Although the terms demonstration or solution are not defined, they are clearly distinguished hence drawing explicitly attention of the reader on their different functions. Other subsections are corollary and remarks (or scholia), the former has the structure of a remark with no specific subdivisions although it contains two different parts, the statement of the corollary and its justification.

---

[10] « […] mon but a été de faire des éléments très rigoureux. J'ai suivi d'assez près la méthode des éléments d'Euclide, et celle du livre d'Archimède de Sphaera et Cylindro : mais en tâchant d'égaler ou même de surpasser mes modèles d'exactitude, j'ai voulu aussi épargner la peine du lecteur autant qu'il m'a été possible, et j'ai fait mes efforts pour donner aux démonstrations toute la clarté et la brièveté que le sujet comporte. Je suppose que le lecteur ait connaissance de la théorie des proportions, qu'on trouve expliquée dans les traités ordinaires d'arithmétique ou d'algèbre ; je lui suppose même la connaissance des premières règles de l'algèbre […] Pour nous, qui avons cet instrument de plus qu'eux, nous aurions tort de n'en pas faire usage s'il en peut résulter une plus grande efficacité. […] et il serait puéril d'employer toujours une méthode laborieuse tandis qu'on peut lui en substituer une beaucoup plus simple et aussi sûre. » (Legendre, 1794, p. vj-viij)

[11] « La géométrie est peut-être, de toutes les parties des mathématiques, celle que l'on doit apprendre la première ; elle me parait très propre à intéresser les enfants, pourvu qu'on le leur présente principalement par rapport à ses applications, soit sur le papier, soit sur le terrain. » (1799/1804 p.xxix)

[12] Following Reviel Netz (1999, p. 98), Greek mathematicians did not define second-order terms, the metalanguage taking its terms from natural language. This innovation is not a creation of the mathematicians of the 18[th] century, second order terms were defined in Geometry of the 17[th] century, may be before.

Adrien Marie Legendre wrote his Elements during the French Revolution Terror. He was not teaching at the time and took the opportunity of this tragic break to revisit Euclid Elements. It is the work of a mathematician, with a mathematics agenda[13]. Nevertheless, his writing was driven by a concern for simplicity (Barbin, 2007) with in mind an educated public[14]. An important aspect of the book are the Notes which discuss and analyse certain proofs[15] or conceptual difficulties[16]. As a textbook, it is remarkable that the Elements focused as much on Geometry as on proofs, which in themselves have to be understood and learned.

The style and organisation of the text of Lacroix's Elements share the same characteristics which put the issue of proving on the fore. The content is organised in relation to the need for rigorous proofs without unnecessary complexity.  These Lacroix's Elements are a textbook written by a mathematician for an advanced level of education at the time. Aimed at teaching Geometry, it includes an unusually lengthy preface dedicated to the method in mathematics, emphasizing and discussing the logic and rigour of the mathematical discourse. It makes Geometry as much as proof the object of learning, as it is witnessed by this excerpt of the *Preliminary Discourse*:

> « *I would add that one must not neglect to present in the geometric demonstrations, an example of the various forms of reasoning, to show how the rules of Descartes and Pascal are observed, and how the certainty of Geometry results from the precise determination of the objects it considers, and therefore each one which can only be envisaged under a very limited number of faces, lends itself to complete enumerations, which leave no doubt as to the result of the reasoning. Elements of Geometry treated in this way would in some way become excellent elements of logic, and would perhaps be the only ones that should be studied* ». (Lacroix, 1799, p. xxviij).

The criticisms of Euclid's Elements from a mathematical point view and the point of view of their usefulness (i.e. Geometry versus Practical geometry) have accompanied their dissemination from the outset. What changes at the end of the 18th century is the social and political context of the teaching of mathematics. Until then, Geometry was taught to a privileged class of people, mostly adults. That changed in France, after the 1789 Revolution. The teaching of mathematics became part of a national educational policy. It was the case as well in most European countries (E. Barbin & Menghini, 2014, p. 475; Schubring, 2015, p. 242-244).

The need for better educated citizens and workers led nations to organize public educational systems, to develop primary and secondary education, and to establish institutions to train teachers; the "normal schools" as they were called after the name of their German precursors of the 18th century. In France, the first one, known as the *Ecole Normale de l'an III*, was created in 1795, for four months, with the assignment *to teach the art of teaching* to educators whose mission was to create afterwards Normal School in townships[17]. Dedicated to primary education, which had a political priority, these schools trained educators to teach practical geometry. The structuration of education beyond primary education, during the first third of the 19th century, distinguished different teaching of mathematics depending on the political views on the future of students. As Hélène Gispert (2014) commented it: "to each social class 'its' mathematics: formation of the mind versus training for practice" (ibid. section 2.1). The need to legislate on educational contents elicited the epistemological and educational rupture between Geometry and Practical geometry. The theoretical nature of the former and the emphasis on the role of proof – as the expression of deductive reasoning – is without doubt source of this rupture. The resources of the educational system, limiting the duration of studies and subject to society priorities, requested making a choice.

---

[13] The report to the *Conseil des Cinq-cents on the elementary books* submitted to the competition opened by the law of 9 Pluviôse, Year II (Lakanal, 1795), cites the works whose manuscripts were presented and retained, giving the main reasons for the choices. For mathematics, it is noted that works that are « not very rigorous, and not very suitable for accustoming children's minds to exact reasoning » were excluded. The list of works chosen is followed by a note recommending Legendre, which, among the printed books, "must be placed in first place" because "his reputation is not disputed, even by envy". This is the only reason; in my opinion Legendre proposed his Elements without considering the very purpose of the competition.

[14] See (Legendre, 1794, preface)

[15] The famous note on the proof of the invariance of the sum of the angles of a triangle (e.g. Note 2)

[16] See e.g. note XII which starts from a discussion of the definition of equality or similarity between polyhedrons.

The massive development of basic education raised the need for elementary textbooks, hence the search for an efficient didactic transposition considering the students and privileging understanding over convincing (Barbin, 2021, p. 106-107). But the notion of elementariness differs whether one considers primary[18] or secondary education. The former targeted basic literacy and knowledge of practical value for citizens and an industry which was rapidly growing[19], the latter targeted the acquisition of the foundations by students meant to enter higher education. We must have in mind that only a small number of students, mostly boys[20], entered secondary education and were exposed to the teaching of Geometry. However, this teaching had little room compared to the study of Latin and humanities (Gispert, 2002, p. 4). I focus for the rest of this subsection on secondary education.

The writing of Geometry textbooks was driven by *arguments of simplicity*[21] with a focus which varied from authors to authors, addressing the ordering of the theorems (logical versus natural), or the nature of geometrical objects (simple and elementary objects versus compound objects), or the nature of the first ideas (axioms, definitions) or eventually principles and proofs (E. Barbin, 2007, p. 226-236). In this context Legendre's Elements and Lacroix's Elements deserve a special attention, because they were largely used and disseminated until the end of the 19th century. Lacroix's Elements, thanks to the centralized ruling of the educational system, took a monopolistic position in France, while Legendre was the most disseminated internationally being in some countries seen as true rival of Euclid (Schubring, 2007). For instance, Nathalie Sinclair points "an invasion of French mathematics: The geometry textbook of Adrien-Marie Legendre – and textbooks were the defining curriculum then – began taking the place of Euclid at the American universities, and the influence of the British waned." (Sinclair, 2006, p. 17).

With very different arguments developed in their preface, Legendre and Lacroix textbooks evidence the emergence of *mathematical proof – demonstration* in the text – as an essential part of the learning of Geometry in secondary education. This level of education being mostly attended by students ambitioning long studies and the entrance in universities. The tension between the utilitarian and the theoretical nature of the curriculum at this level was not too high; the teaching of geometry had a theoretical coloration. More important was the tension between proof that explain and proof that convince, to use the terms of contemporary discussions (Hanna, 2000), or, with words of the time, the tension between *analysis* and *synthesis*. The « Preliminary comments » of Sylvestre-François Lacroix was dedicated to this issue. This will be one of the issues to be addressed by the 20$^{th}$ century with the democratization of secondary education.

---

[17] (« 10. 30 octobre 1794 (9 brumaire an III). Décret relatif à l'établissement des écoles normales », 1992)

[18] Compulsory school lasted until the age of 12 at that time (Gispert, 2002, p. 6).

[19] The geometry taught in primary schools favoured a concrete approach and know-how by mobilising drawing and manual work (D'Enfert, 2003, p. 7).

[20] In France, during almost all the 19th century, « the lycée "was reserved for 3 to 4% of an age group and for young people only, a paid education, of culture, for which mathematics was disqualified, relegated to the final classes as an element of specialisation. » (Gispert, 2002, p. 4). Most girls were left out secondary education, when they got access at the end of the century it was with a special "female" course of study which could be qualified "second-rate education" (ibid.).

[21] The word simplicity translates here the French *élémentation*—which is of the same family as the word *élément*—introduced by Evelyne Barbin (2021). It refers to « the ordering of a science, here elementary geometry, which seems to be the most appropriate for its teaching. The term "elements" is present in the title of the oldest mathematical work that has come down to us, that of Euclid dated from the 3rd century BC. »  (ibid. p.99).

Patricio Herbst calls the *Era of Text*[22] that period in which: "The study of geometry was done through reading and reproducing a text; such work would train the reasoning faculties of students. But, the texts do not hint at the existence of official mechanisms to verify or steer the evolution of students' reasoning." (P. G. Herbst, 2002b, p. 288 ff.). This discrepancy between intentions and means received particular attention in the USA in the last decade of the 19th century. In the USA, education has never been driven by a single state institution, but by many local agencies with a wide range of organisational, pedagogical and epistemological views. This created problems for all disciplines in recruiting students at college level, which the National Education Association addressed by creating the Committee of Ten in 1892. The committee's role was to help school districts and private schools make changes by providing arguments to support decisions based on what universities would expect. In mathematics (P. G. Herbst, 2002b, sect. 2), the diversity of approaches highlighted the tension between educating the mind and transmitting knowledge.

The Mathematics Conference, convened by the Committee of Ten in 1893, reached a consensus that the education of the mind of secondary school students should take priority, and that Geometry should have a role in the development of reasoning skills. "[It] recommended changes to the geometry curriculum to accommodate the tension between training mental faculties [i.e. justification] and imparting culturally valued geometric knowledge" (P. G. Herbst, 2002b, p. 295). How to teach students how to construct proofs in mathematics became an explicit question. The pioneering work of George Albert Wentworth had provided an answer by proposing a norm for layout in which "each distinct assertion in the demonstration, and each particular direction in the constructions of the figure, begins a new line; and in no case is it necessary to turn the page in reading a demonstration" (Wentworth, 1877, p. iv). The preface to the third edition[23] includes a section "For the teacher" with, among other recommendations:

> "The teacher is likewise advised to give frequent written examinations. These should not be too difficult, and sufficient time be allowed for accurately constructing the figures, for choosing the best language and for determining the best arrangement." (ibid. p.vi).

This was a precursor to the two-column form that dominated geometry teaching in the USA in the 20th century.

However, for the time being, we can notice that at the end of the 19$^{th}$ century, *proof was* a *named but implicit content to be taught* while teaching Geometry remained the explicit agenda.

## The first part of the 20$^{th}$ century, proof and the formation of scientific mind

The rapid development of the industrial economy and of manufacturing engineering in the early 20$^{th}$ century highlighted the need to improve mathematical literacy and skills of the workforce at all levels. This concern is international (Nabonnand, 2007). The turn of the 20$^{th}$ century is also the time of the international organization of mathematicians with establishment of "The International

---

[22] Patricio Herbst coined this expression for the USA, in my opinion it can be extended to Europe.
[23] The third edition (1881) of this text book is available at [https://hdl.handle.net/2027/hvd.32044097014377].

Commission on the Teaching of Mathematics" (IMUK), in 1908[24] and the creation of the International Mathematical Union (IMU), in 1920. The creation of IMUK demonstrates the international concern for the development of the teaching of mathematics. One of the first decisions of this Commission is « to survey and publish a general report on current trends in mathematics education in the various countries » (H. F. Fehr, 1908, p. 8). It requires the survey to consider applied as well as pure mathematics, and recommends that it focuses on principles which should inspire the teacher, but it leaves aside curricula which are the responsibility of nations.

The question of rigor received a special attention. A report on this issue is presented to the IMUK delegates. The rapporteur, Guido Castelnuovo, proposed to limit the discussion to the upper secondary schools and to the teaching of geometry. The topic is the extent to which the systematic presentation of mathematics can be considered. A classification of the degree of rigor is proposed:

> « A) *Entirely logical method* – All axioms are stated; their independence is discussed; further development is rigorously logical. No appeal is made to intuition; primitive notions (point, etc.) are subject only to the condition of satisfying the axioms.
>
> B) *Empirical foundations, logical development* – From the observation of real space, we deduce the primitive propositions on which the following logical development is based. The following logical development is based on these propositions -- three subgroups can be distinguished: $B_A$ all axioms are stated, $B_B$ some of the axioms are stated, $B_C$ only those axioms which are not absolutely self-evident are stated.
>
> C) *Intuitive considerations alternate with the deductive method* – Evidence is used whenever appropriate, without it being clear what is admitted and what is demonstrated.
>
> D) *Intuitive-experimental method* – Theorems are presented as facts that are intuitive or can be demonstrated by experience, without the logical connection between these facts being apparent. »
> (H. Fehr, 1911, p. 462)

It appears that no country chose A or D. Guido Castelnuovo noticed that Roman nations and the UK prefer B, and that German nations are closer to C. His comments suggested an influence of culture and possibly of the economical context (esp. « industrialism »).

The exchanges underline the importance of a *non-excessive rigour* considering the average « intelligence » of the students: Rigour must be compatible with teaching, and if learning geometry favours the development of logical reasoning, it is not necessary to go as far as the installation of an axiomatic (ibid. pp. 465-466). The reference to the psychology of the young if frequent in the justification of the choices made by nations. It is proposed to discuss in the future the organisation of the teaching of geometry and to study its psychological grounds.

From now on modern mathematical education became a national stake in most nations, some henceforth searching for curricula balancing the applicative and theoretical value of the teaching of Geometry (e.g. González & Herbst, 2006).

The classification of the arguments of the Geometry course by Gloriana González and Patricio Herbst's facilitates distinguishing and understanding the different rationales that shape the didactic transposition of mathematical proof (ibid. p. 13):

---

[24] In 1954 IMUK changed its name for ICMI (International Commission on Mathematical Instruction) with the mission of "[the] conduct of the activities of IMU, bearing on mathematical and scientific education" (Furinghetti & Giacardi, 2008).

1. a *formal argument* that defines the study of geometry as a case of learning logical reasoning through the practice and application of deduction;
2. a *utilitarian argument* that geometry would provide tools for the future work or non-mathematical studies;
3. a *mathematical argument* justifies the study of geometry as an opportunity to experience the work of doing mathematics;
4. an *intuitive argument* aligns the geometry course with opportunities to learn a language that would allow students to model the world;

These arguments are not mutually exclusive, several of them could contribute to the didactic transposition of geometry, but with different weights. With regards to mathematical proof, the "formal argument" and the "mathematical argument" support the *raison d'être* of its teaching. The other two arguments are less decisive because its teaching cannot claim to provide a model of proving for all areas of knowledge, whether scientific or practical as the fourth argument suggests. It happens that proof is often seen as an obstacle to the learning of geometry for its lack of practical value. On the contrary, the arguments of the mathematicians of the beginning of the century is that mathematical proof is constitutive of Geometry as a paragon of mathematics as science. This concern for the teaching of mathematics as a *science* is well illustrated by a comment of Giuseppe Veronese after the Castelnuovo reporting: « If industrialism or general utilitarianism were indeed a dominant influence in middle school education, mathematicians would have to fight it. » (H. Fehr, 1911, p. 465).

At the turn of the 20th century, it was clear that the learning of deductive reasoning is an important educational objective. Mathematics got more importance than the Humanities which were until then the educational priority. It was meant to play a privileged role in the formation of *the scientific mind*, as the French 1946 General instructions called it. Geometry was the elected domain: « *It is important to make the difference felt very early on between the certainty given by the geometric method and that resulting from the experimental method: it is on this condition that the need will develop for demonstration* » (French Instructions of 2 September 1925). However, achieving this objective for early grades proved to be a challenge:

> « *But is it possible to ensure understanding of mathematics among young pupils, especially in the sixth and fifth grades? The question is still being discussed and the instructions that followed the 1902 reform went so far as to prohibit theoretical explanations on certain points and in certain classes. One had to be content to have the rules learned and applied, for a well-fixed mechanism.* » (French Instructions, 1923).

Decision makers searched for solutions introducing pedagogical recommendations, as for example le following:

> « *As the hypotheses or data are recorded on the figure itself, by the means most likely to ensure their immediate vision and scope, the teacher would slowly deduce, with the help of the class if possible, the hypotheses or data; he or she would summarise the results acquired at each moment and have the pupils formulate them themselves. The pupils would no longer be confused by the assembly of terms accumulated in synthetic statements whose formation would be partly their own work. They would stop more at the most important ones: the theorems would take shape at the right moment; they would be fixed in the memory by the usual procedures.* » (French Instruction of September 2, 1925).

The limitations of such an approach were anticipated: « faith in the correctness of the rule and confidence in the authority of the teacher contribute to delaying the awakening of the critical sense » (ibid.).

The search for the most efficient way of teaching mathematical proof was constant along this first part of the 20th century. The driving idea was to engage students in problem solving and managing a seamless transition from the manipulation of objects to reasoning on abstract representations. We may say that decision makers and mathematicians understood the rupture between practical thinking and theoretical thinking, but looked for a way to bypass it instead of facing it. Here is another evidence of this approach:

> « Guided by the teacher and first carrying out concrete operations applied to given objects, the child will acquire the abstract notion of an operation of a well-defined nature but concerning an indeterminate element. Then he will become capable of imagining that he is applying another operation to the result of the first one without having carried it out. Finally, designing the continuation of the mechanisms of the operations thus defined, he will be able to predict certain properties of the results: he will have carried out his first demonstration. » (French decree of 20 July 1960).

During this first part of the 20th century, the teaching of Geometry included exercises and problems, providing students possibilities to craft proofs either to achieve simple deductive tasks of one or two steps, or more complex ones requiring students to engage in problem-solving; however, these more complex problems were often cut into parts making them easier to solve. Eventually, although learning was supported by more significant activities, the basic approach consisted of observing, reading and replicating proofs.

Patricio Herbst identifies the results of the Committee of 10 as a turning point between pedagogical approaches, following which students had the opportunity to use their reasoning on corollaries of theorems or theorems not to be included in the main text of the course, namely the Era of Originals (P. G. Herbst, 2002b, sect. 4), and approaches, which proposed activities aimed at learning what a proof is and at practising proving, namely the Era of Exercise (ibid. sect. 6). This evolution was accompanied by the development of a distinctive didactic tool, of which the layout standard proposed by Wentworth was a precursor: the two-column proof. The pattern of this layout is made of two lines forming a T. Above the horizontal line is written the statement to be proved, below that line, separated by a vertical line, two columns display the proof writing with on the right-hand side the sequence of inferences and on the left-hand side the warrant of each of them. It is commonly acknowledged that distinctive layout was first used in the second edition of the Schultze and Sevenoak' Geometry textbook in 1913[25]. It was meant to be a tool for the students as well as for the teacher:

> "[This arrangement in two columns] seems to emphasize more strongly the necessity of giving a reason for each statement made, and it saves time when the teacher is inspecting and correcting written work." (Shibli (1932) comment quoted by P. G. Herbst, 2002b, p. 297)

Two-column proofs brought stability to the Geometry course in the USA, but over emphasizing a formal display of the logical structure of proofs, it tended to hide its role in knowledge construction. (P. G. Herbst, 2002b, sect. 7).

---

[25] (for a review in support of this book see O'Reilly, 1902)

During this first half of the 20th century, *proof got the explicit status of a mathematical tool to be taught and learned* but which learning was induced by the learning of Geometry which had an exclusive focus in curricula. Following the Patricio Herbst's formula: "To know geometry and to be able to prove the theorems of geometry were indistinguishable." (P. G. Herbst, 2002b, p. 289).

## The second part of the 20th century, proof liberated from Geometry

In the middle of the 20th century, mathematics was present in all domains from natural to human and social sciences. Mathematical competences imposed themselves for their key role in the development of modern industry and economic sectors (D'Enfert & Gispert, 2011, p. 30; Gispert, 2002, p. 9). Mathematics was emerging as a universal language for accessing knowledge. Again, countries expressed the need for people better educated in mathematics.

On the academic side, the increasing distance between school mathematics and mathematics as a science and the intellectual influence of the French mathematicians' group Bourbaki and of its Elements, not to mention the Sputnik crisis in the US, led to a definitive rupture with the text of Geometry inherited from Euclid.

The Royaumont seminar on school mathematics in December 1959[26] gave the direction for the future. The Euclidean text was definitely considered obsolete from both a mathematical and a pedagogical perspective, but this left the mathematics education community with more problems than solutions. Dieudonné exclamation *à bas Euclide* attracted the attention of the general public, but didn't account for the discussions on what the desired evolution should be like. There was a large consensus on the final goal of geometry instruction, viz. that after the early stages of intuitive learning, there should come « the breaking of the bridge with reality – that is, the development of an abstract theory. » (OECD[27] quoted by (Bock & Vanpaemel, 2015, p. 159)). The OEEC[28] official report was published in 1961, under the title "New Thinking in School Mathematics" (popularised as "New Math"), and "Mathématiques nouvelles" for the French version (ibid. p. 163).

The movement spurred on by some mathematicians looked for an epistemological break. It led to two guiding principles for the design of new curricula (D'Enfert & Gispert, 2011, p. 35):

1. mathematics is a deductive science, not an experimental science.
2. mathematics forms a theory[29] which must bring together under the same structure knowledge that has hitherto been presented in a scattered manner.

Geometry had to evolve. This evolution was radical in the US where Geometry was relegated to teacher training and became optional in schools (Sinclair, 2006, p. 73-74). In other places it was deeply transformed. In France it maintained itself in high school curricula but with a new face in which the study of geometrical objects gave way to that of structures, much attention being paid not to confuse the concrete world with its mathematical model using different terminology. The main influences leading to these evolutions were of different origins in the different countries. Dirk De Bock and Geert Vanpaemel analysing the OEEC seminar at Royaumond noticed that "In France the demand for modernization came from the universities and was aimed at introducing modern 'Bourbaki mathematics' in secondary schools; in the United States the renewal of mathematics education was urged by industry and politics, and aimed at the modernization of teaching methods" (Bock & Vanpaemel, 2015, p. 167). Moreover, the participants to this seminar which was decisive for the mid-20[th] century evolution were not all aligned on the Dieudonné's position often quoted as the slogan of Modern mathematics[30]. A much more balanced approach to the reform was being proposed (ibid. 152). Dieudonné himself outlined a curriculum "quite concrete, roughly starting from 'experimental' mathematics, concentrating on techniques and practical work, to a rigorous, axiomatic treatment of two- and three-dimensional space." (ibid. 157). This has found a direct translation in the French official documentation of curricula: « Success will be achieved when the pupil, having become aware of the difference between an experimental verification, even if it is repeated a hundred times, and a demonstration, comes not to be satisfied with the first and to demand the second. » (Decree of 20 July 1960). However, the transition from the so-called experimental verification to mathematical proof (i.e. demonstration) was radical; classically, the *rupture* happened at the 8[th] grade. As emphasizes Gert Schubring (2015), France engaged a reform with no consideration for the needs of the different students' orientation towards vocational studies or university studies: "[For the 8th grade and the 9th grade], common to all the diverse curricular directions, the Commission had planned to teach the same contents, and according to the same spirit and methodology – conceiving this exclusively from the logic of a curriculum for those who would continue to university studies." (ibid. p. 250).

---

[26] "The origin of the Seminar must perhaps be placed as far back as the 1952 meeting of CIEAEM in La Rochette par Melun on "mathematical and mental structures", which had brought together Dieudonné, Choquet and Servais in dialogue with psychologist Jean Piaget and philosopher Ferdinand Gonseth. In several countries the reform of school mathematics was well underway by 1959, with a large number of specialist meetings on a regular basis" (Bock & Vanpaemel, 2015, p. 165-166).

[27] Organisation for Economic Co-operation and Development

[28] Organisation for European Economic Co-operation

[29] In French, even the name of mathematics changes, losing the mark of the plural to become "La Mathématique", on the model of "la physique" which is a singular word in French (but a plural in English "physics"). This did not last for long.

Axiomatic, which I consider the true heritage from Euclid, is of a paramount importance in the reform. Proof is its backbone, axioms and definitions are its ground. Geometry is the privileged terrain for being acculturated to this new epistemology: « born from experience, [it] should appear to students as a true mathematical theory » at the end of the 8[th] year (French Decree of 22 July 1971). However, this does not apply to geometry only but to all mathematics. The learning of mathematics as a discipline should train students in *deductive thinking*, encourage them to be *rigorous in logic*, teach them *to build a chain of deductions*, to develop – in a constructive way – their critical mind[31].

The theoretical orientation of the New Math was criticised internationally, by both mathematics educator and policy makers, with – at least – two arguments: the break between mathematics and its applications including its use by other scientific disciplines, and its irrelevance for a large population of students. This was well expressed by this judgement from French the pioneer of applied mathematics, Jean Kuntzmann (1976, p. 157):

> « One could regret that it is not possible to lead the students to the deductive phase. In reality, nothing is lost because the philosophy of this phase: autonomy of mathematics organizing itself in view of its own objectives, is perfectly useless to the average Frenchman[32] (even of the year 2000). I affirm very clearly that for the average Frenchman, therefore for the teaching of the first cycle, the suitable philosophy is that of the conceptual phase. That is to say:
>
> > - Duality situation-model, fundamental for the uses of the mathematics;
> > - Training in logical reasoning but without going as far as organized deductive theory (one will meet in everyday life occasions to reason, but few constituted deductive theory). »

---

[30] The section of the seminar report subtitled *Sharp Controversy Provoked*. "After some discussion, both groups modified their positions on the programme and reached general agreement on a set of proposals which did not remove Euclid entirely from the secondary-school curriculum" (OEEC, 1961, p. 47) (quoted by (Bock & Vanpaemel, 2015, p. 157)).

[31] Excerpt from the French 29 April 1977 circular.

[32] Only a third of the students engage then in the secondary school (D'Enfert & Gispert, 2011, p. 40).

In the USA, the New Math movement declined in the early 1970s under the pressure of the public concern reflected by the catch phrase "Back to basics", a decrease of students' performances and the criticisms of some leading mathematicians (Sinclair, 2006, p. 108 sqq). In fact, it had strong opponents since the very beginning (e.g. Goodstein, 1962). Moreover, the rapid development of computer-based technologies and the growing evidence that computers will change the role of and demand on mathematics called for questioning curricula. The National Advisory Committee on Mathematics Education (NACOME) published in 1975 recommendations for an evolution of compulsory school curricula which did not reject all the contributions of the New Math but redefined the "basics" (Hill, 1976) putting on the fore applications of mathematics, statistics and probability, the use of calculators and computers. Remarkably, Geometry was not part of the Basics mathematical skills listed by the NACOME report but the work for reforms following the recommendations included it (Sinclair, 2006, p. 111). Proof and logical reasoning were down played in favour of a "wider conception of geometry" giving room to visualisation and intuition (ibid. p.113); a move echoing the call of leading opponents for abandoning the "chief innovation" of the New Math: the logical approach to the teaching of mathematics (Kline, 1976, p. 451-453). By the end of the century, the NCTM Principles and Standards for School Mathematics (PSSM) made a synthesis of both positions:

> "Geometry has long been regarded as the place in high school where students learn to prove geometric theorems. The Geometry Standard takes a broader view of the power of geometry by calling on students to analyse characteristics of geometric shapes and make mathematical arguments about the geometric relationship, as well as to use visualization, spatial reasoning, and geometric modelling to solve problems. Geometry is a natural area of mathematics for the development of students' reasoning and justification skills." (NCTM, 2000, p. 3)

The New Math movement faded away in all countries by the early 1980s. This end was a consequence of the constant criticism of curricula which, burrowing the words of José Manuel Matos[33], "render mathematics hermetic" either for students or their parents and most stakeholders as well, and a consequence of the lack of consensus among decision makers, mathematics educators and teachers as well as mathematicians themselves. It happens that the gap between the reform and the reality of the classrooms was such that it was not surprising to see the New Math cohabiting with the preceding classical teaching, or even being ignored, in particular where the educational system left enough autonomy to schools and teachers. In France, where mathematicians had a peculiar responsibility in launching the movement, the tension between the protagonists led to the creation of a *Union of the users of mathematics* opposing the reform[34].

The following reforms did not return back to the old curricula. The conception of mathematics, its scholastic organisation, the content of its various components and the focus on the *mathematical activity* rather than on the *text of mathematics* moved significantly. The debate initiated in the USA in the beginning of the 1970s lasted ten years before the release of consensual recommendations in the form of *An Agenda for Action* endorsed by the National Council of Teachers of Mathematics in 1980:

> "*An Agenda for Action* (NCTM 1980), recommended that problem solving should be the focus of school mathematics, that basic skills should be defined more broadly than simple arithmetic and algebraic calculation, that calculators and computers should be used at all

---

[33] In (Ausejo & Matos, 2014, p. 298)

[34] UPUM created in 1972 has disappeared since then with the New math.

grade levels, and that more mathematics should be required of students." (Fey and Graeber 2003, p. 553 quoted by (Kilpatrick, 2014, p. 331))

This NCTM Agenda led to a questioning of the teaching of proof which, through a metacognitive shift[35], had in practice evolved into the teaching of the two-column proof technique. It eventually motivated the recommendation in the 1989 Standards to increase attention to "deductive arguments expressed orally and in sentence or paragraph form" (p. 126) and to decrease attention to "two-column proofs" (p. 127)." (Quoted by P. Herbst, 1999).

The creation of the *International Congress on Mathematical Education* in 1969 and of the *international group for the Psychology of Mathematics Education* in 1976 favoured the dissemination of the various positions and ideas, allowing international debates and exchanges between mathematics educators and teachers at an international level. Thus, the post-New Math orientation of the curricula was rather similar in most of the countries. For example, in South Asia:

> "The Math Reforms lasted for 12 years, ending in the early 1980s, when it was realized they did not work and had to be stopped. Although many new topics introduced during the Math Reforms stayed on (e.g., Venn diagrams and statistics), the formal approach in teaching mathematics was replaced by the so-called problem-solving approach. In the years that followed, change in content was minor. *The major change was in the teaching approach used in the classroom*." (Lee, 2014, p. 388- my emphasis)

Similar lines have been written in several other educational contexts[36].

In France – which is an interesting case given the constraints of a strong and centralized educational administration and the influential position of professional mathematicians – the 1980s post-New Math era prepared a rupture in the mathematics education policies of the 21st century:

> "Moreover, challenged since the early 1970s, including by its supporters who believed the reform did not correspond to their recommendations, the "modern mathematics" reform was abandoned in the early 1980s in favour of a teaching method that, envisioning mathematics in the diversity of its applications, placed the accent on problem solving and favoured "applied" components of the discipline. These two aspects now occupy a central place in mathematics teaching. At the same time, since the early 2000s, there has developed the ambition to make mathematics into a subject that allows students to throw themselves into a true research program, capable of developing their abilities to reason and argue but also to experiment and imagine." (Gispert, 2014, p. 239)

But, Geometry, at the end of the 20th century, had lost its special and somewhat isolated position in mathematics curricula. It is no longer the Geometry of the Euclidean educational tradition, nor is it the Geometry of the New Math. The study of geometrical figure-objects is part of its educational content along with geometrical transformations, with sometimes important differences among national curricula but this orientation is shared.

The emphasis is on making mathematics meaningful and on understanding, *mathematical proof* loses ground[37] to *deductive reasoning*, which opens up a broader conception of validation in the learning

---

[35] (Brousseau, 1997, p. 26 ff.)

[36] For references, see for example (Karp & Schubring, 2014, part. V)

[37] This claim may have exceptions, as it is the case of Russia where "Russian schoolchildren of the 1980s–1990s, and even the schoolchildren of today, spend much more time than their Western counterparts on algebraic transformations and proofs." (Karp, 2014, p. 320)

and teaching of mathematics. *One observes an epistemological revolution more than a new educational one.*

## The 21st century, proof for all grades

Since 1995, the "Trends in International Mathematics and Science Study" (TIMSS) provides an instrument to get a picture of the institutional views of mathematics education[38]. Its objective is to study 4th and 8th graders' achievements in mathematics and science in all participating countries[39]. Since 2003, the publication of the assessment framework gives a view on teaching and learning which pays attention to including goals considered important in a significant number of countries. It can be seen as a consensual conception of the essential basis of the curriculum, although there are many differences, some of which are substantial. These documents provide a reliable basis for getting an idea of how proof and proving are perceived in the early 21st century.

The design of the TIMSS assessments distinguishes two types of domains, the *content domains* and the *cognitive domains*: "The content domains define the specific mathematics subject matter covered by the assessment, and the cognitive domains define the sets of behaviours expected of students as they engage with the mathematics content." (TIMSS 2003 / O'Connor et al., 2003, p. 9)[40].

Issues related to validation are addressed in the sub-domain "reasoning" which is presented as follows:

> "*Reasoning mathematically* involves the capacity for logical, systematic thinking. It includes intuitive and inductive reasoning based on patterns and regularities than can be used to arrive at solutions to non-routine problems. Non-routine problems are problems that are very likely to be unfamiliar to students. They make cognitive demands over and above those needed for solution of routine problems, even when the knowledge and skills required for their solution have been learned. Non-routine problems may be purely mathematical or may have real-life setting. Both types of items involve transfer of knowledge and skills to new situations, and interactions among reasoning skills are usually a feature.
> Most of the other behaviours listed within the reasoning domain are those that may be drawn on in thinking about and solving such problems, by each by itself represents a valuable outcome of mathematics education, with the potential to influence learners' thinking more generally. For example, reasoning involves the ability to observe and make conjectures. *It also involves making deductions based on specific assumptions and rules, and justifying results*." (ibid. p.26 – my emphasis)

The "cognitive domains" includes "Knowing", "Applying" and "Reasoning". Each of these is characterised by "a list of objectives covered in a majority of the participating countries, at either grade 4 or 8." (ibid. p. 9). In the case of Reasoning, the objectives expressed in behavioural terms are: *Analyse, Generalize, Synthetize/Integrate, Justify*. The latter was labelled *Justify/Prove* in 2003, but only *Justify*[41] remained for the following assessment campaigns.

*Table 1 - Expression of the objective "Justify" in the Assessment framework document of the TIMSS from 2003 to 2019*

| 2003 | Justify/Prove | "Provide evidence for the validity of an action or the truth of a statement by reference to mathematical results or properties; develop mathematical |
|------|---------------|--------------------------------------------------------------------------------------------------------|

---

[38] There are two major international assessment campaigns, PISA and TIMSS. The former assesses the learning achievement of 15-year-olds. The second one does it for 4th and 8th graders. I think the latter is more relevant to the issues I am addressing here.

[39] They were 64 in 2019.

[40] Similarly, the frameworks of the TIMSS 1995 and of the TIMSS 1999 assessments included *content areas* and *performance expectations*.

[41] The expression « making deductions » in the 2003 document, became « making logical deductions" in the 2007 document and remained stable then.

|  |  | arguments to prove or disprove statements, given relevant information." (TIMSS 2003 p. 33) |
| --- | --- | --- |
| 2007 | Justify | "Provide a justification for the truth or the falsity of a statement by reference to mathematical results or properties" (TIMSS 2007 /Mullis et al., 2007, p. 38) |
| 2008 | Justify | "Provide a justification for the truth or falsity of a statement by reference to mathematical results or properties." (TIMSS 2008 / Garden et al., 2008, p. 22) |
| 2011 | Justify | "Provide a justification by reference to known mathematical results or properties." (TIMSS 2011 / Mullis et al., 2009, p. 46) |
| 2015 | Justify | "Provide mathematical arguments to support a strategy or solution." (TIMSS 2015 / Mullis & Martin, 2014, p. 27) |
| 2019 | Justify | "Provide mathematical arguments to support a strategy or solution." (TIMSS 2019 / Mullis et al., 2017, p. 24) |

The 2003 TIMSS assessment framework associates justify and prove, however "to prove" must not be interpreted as providing a *mathematical proof* but as providing *mathematical arguments*. This interpretation is coherent with the need to have a formulation adequate to the 4$^{th}$ grade as well as for the 8$^{th}$ grade. "Proving mathematically" appears explicitly only at the latter grade in most curricula, when it does.

The reference to mathematical proof being abandoned, the keyword which is chosen is "justification" with specific requirements: "reference to mathematical results or properties" (TIMSS 2007, TIMSS 2008), "reference to known mathematical results or properties" (TIMSS 2011). Then comes back the key expression "mathematical argument" (TIMSS 2015, TIMSS 2019) in a short and allusive statement compared with the preceding formulation.

Moreover, the section of the TIMSS 2003 document, entitled "Communicating mathematically" disappears in the following editions. It said that "Communication is fundamental to each of the other categories of knowing facts and procedures, using concepts, solving routine problems, and reasoning, and students' communication in and about mathematics should be regarded as assessable in each of these areas." (ibid p. 34).

These TIMSS assessment framework documents reflect on the one hand an institutional vision of a detachment of the teaching of proof from Geometry, and on the other hand an objective to introduce the learning of a proper way to address the question of "truth" in the mathematics classroom. This evolution of TIMSS witnesses a trend of curricula in a large number of countries. It strengthens a scholastic epistemological revolution which the following quote clearly exemplifies:

> "One hallmark of mathematical understanding is the ability to justify, *in a way appropriate to the student's mathematical maturity*, why a particular mathematical statement is true or where a mathematical rule comes from." (NGA Center & CCSSO, 2010, p. 4 - my emphasis).

It should be noted that at that time research in mathematics education was reaching academic maturity. The scientific community had acquired the necessary professional tools to consolidate and to establish a science that had been asserting itself since the mid-1970s: international journals and conferences with clear scientific policies and quality control. The ICME conferences and the ICMI studies, under the umbrella of IMU, maintained the links between researchers in mathematics education and mathematicians. The creation of the ICMI awards in 2000 in keeping with the tradition of the mathematical community is another evidence. All these are indicators that researchers in mathematics education are now members of the stakeholders' community[42] which contributes to the didactic transposition of the mathematics to be taught and learned.

Researchers in mathematics education have fully embraced the problems of teaching proof, which they claim essential for the learning of mathematics. The number of articles and conference papers on the learning and teaching of proof and proving in the mathematics classroom has impressively increased since the pioneer work of Alan Bell (1976). The number of working groups and of study groups speaks for the dynamic of the community, including edited books contributing to shorten gaps in research (e.g. Boero, 2007; Hanna & de Villiers, 2012; Reid & Knipping, 2010).

In her preface to the book "Theorems in School" (Boero, 2007), Gila Hanna writes in clear words that "[…] proof deserves a prominent place in the curriculum because it continues to be a central feature of mathematics itself, as the preferred method of verification, and because it is a valuable tool for promoting mathematics understanding" (ibid.p.3). Paolo Boero idea of this book was born in the context of the 21st PME conference (Pehkonen, 1997, vol. I-p. 179-198) following a research forum which demonstrated "the renewed interest for proof and proving in mathematics education" and that of "the reconsideration of the importance of proof in mathematics education was leading to important changes in the orientation for the curricula in different all over the world" (Boero, 2007, p. 20).

In 2007, ICMI launched its 19th study on "Proof and proving in mathematics education" (Hanna & de Villiers, 2012/2021) which *Discussion document* introduces the idea of "developmental proof" (ibid. p. 444):

> "The study will consider the role of proof and proving in mathematics education, in part as a precursor for disciplinary proof (in its various forms) as used by mathematicians but mainly in terms of developmental proof, which grows in sophistication as the learner matures towards coherent conceptions. Sometimes the development involves building on the learners' perceptions and actions in order to increase their sophistication. Sometimes it builds on the learners' use of arithmetic or algebraic symbols to calculate and manipulate symbolism in order to deduce consequences. To formulate and communicate these ideas require a simultaneous development of sophistication in action, perception and language.
>
> The study's conception of "developmental proof" has three major features:
>
> 1. Proof and proving in school curricula have the potential to provide a long-term link with the discipline of proof shared by mathematicians.
> 2. Proof and proving can provide a way of thinking that deepens mathematical understanding and the broader nature of human reasoning.
> 3. Proof and proving are at once foundational and complex, and should be gradually developed starting in the early grades."

These features resonate with the modal arguments – for the Geometry course – proposed by Gloriana Gonzáles and Patricio Herbst (2006) namely the *formal*, *utilitarian*, *mathematical*. This categorization can be reused substituting *proof and proving* for *geometry* without losing its relevance.

The analysis from a developmental proof perspective of contemporary institutional texts that provide teachers with comments and pedagogical indications shows the necessary modulation of this teaching considering both the school levels, societal needs and its possibility in relation to mathematical requirements. Let us take the case of France (Balacheff, 2022): from grades 1 to 3, the students' discourse must be argued and based on observations and research and not on beliefs; from grades 4 to 6, the teaching must contribute to students building the idea of proof and argumentation (e.g. by moving gradually from empirical validations to validations based solely on reasoning and argumentation). In grades 7-9, the challenge is to move from inductive to deductive reasoning, and to put this deductive reasoning into the form of a communicable proof (i.e. a *demonstration* in the French text)[43]. This is, with variations, what is observed internationally and reflected in TIMSS. The main point of divergence is the point at which acculturation to the socio-mathematical norm of mathematical proof is targeted; in many countries this is left to the upper secondary school level and often, but not always, to the learning of Geometry (e.g. in the US, Jones & Herbst, 2012, p. 263).

---

[42] The noosphere is "the sphere of those who "think" about teaching, an intermediary between the teaching system and society." (Chevallard & Bosch, 2014, p. MS 1)

Although the institutions stress that the teaching and learning of proof should not be confined to Geometry, this domain remains an ecological niche for achieving this goal. Writing a "Spotlight on the Standard" for the NCTM journal "Teaching Mathematics in the Middle School", Edward A. Silver (2000) made the remark that "Although many middle school students love to argue (about almost anything!), they need to learn to argue effectively in mathematics. The study of geometry offers many opportunities to gain experience with mathematical argumentation and proof" (ibid. p.23). As it happens, the discussion document of the 19th ICMI study reserves a special place for geometry when questioning the research community on the relation between proof and empirical science, "given geometry deals with empirical statements about the surrounding space as well as with a theoretical system about space" (Hanna & de Villiers, 2012, p. 451). The availability of Dynamic geometry microworlds resonates with this questioning, at a point where the designers of the study devote to it a whole section under the title "Dynamic geometry software and transition to proof" which first question is: "To what extent can explorations within DGS foster a transition to the formal aspects of proof? What kinds of didactic engineering can trigger and enhance such support? What specific actions by students could support this transition?" (ibid. p. 449)

## An epistemological rupture in need for an instructional[44] bridge

Coming back over several decades of study of the problems raised by the teaching of proof in the mathematics classroom, I realised that I lacked a comprehensive study of the history of this instructional objective. Research on the history of mathematics education provides us with a number of information and analysis, especially in the literature on the history of the teaching of geometry. My objective was to gather this information, adding when needed some evidence from primary resources, and to structure it in order to get a picture of the history of the teaching of proof which could be useful to carry out a study which is still to be done.

To conclude these notes, I will first outline what can be retained from a historical point of view, and then present arguments in favour of searching for a precursor concept of mathematical proof that allows the question of truth to be addressed in the early teaching of mathematics.

## Milestones in a long journey in search of a solution for the early learning of proof

From c. 300 BC to the late 18th-century, Euclid's Elements stood as the model of the text of "scholarly knowledge"[45] of Geometry, which is the material for the process of the didactic transposition. At the turn of the 19th century, under the growing criticisms of the Elements and a development of the progressive ideas on education, started a scientific work pursuing the objective of writing a rigorous presentation of Geometry but with the associated intention to facilitating its understanding for a reader eager to learn it. Yet, it is a transpositive work insofar as it "*improves* the organization of knowledge and makes it more understandable, structured and accurate, to the point that the knowledge originally transposed is itself bettered." (Chevallard & Bosch, 2014, p. MS 2).

---

[43] It is interesting to notice that the institutional discourse avoids the reference to abductive reasoning – Polya would write *plausible reasoning*. Possibly in fear of opening room to severe logical errors (often induced by natural reasoning), although abductive reasoning has a heuristic value and is the source of creative ideas.

[44] Following and understanding the remarks of Keith Jones and Patricio Herbst (2012, p. 261-262), I use *instructional* as an interpretation – if not a translation – of the French word *didactique* although it has a larger denotation; but instructional and didactical objective being tightly related this does not open serious misunderstanding in the context of this Note.

[45] (Chevallard & Bosch, 2014, fig. 1)

Until the end of the 18$^{th}$ century, Euclid's geometry was taught to a privileged class of the society, in particular – burrowing the words of Legendre – to those who were devoted to mathematics; they were mainly adults. Along the 19$^{th}$ century, with the development of national policies, the challenge of teaching Geometry initiated a transposition process which progressively included learning issues with concerns for ever wider segments of the population. The cases of Dechalles, Legendre, Clairaut and Lacroix are significant milestones of this early period of the history of the teaching of Geometry and the way the issue of Proof was discussed and addressed. The former and the latter wrote textbooks with an explicit critical position towards the seminal Elements, while the other two were first writing treatises. But the work of the four of them evidences the awareness of the epistemic complexity of the project:

- the tension between proving and explaining
- the conflict between the abstractness of Geometry as a theoretical construction and its practical value as a tool for numerous human activities.

Proof and proving were a core concern for their role in the understanding of Geometry and in establishing the truth of geometrical statements. Legendre and Lacroix included demonstration among the metamathematical terms which understanding was necessary for the presentation of their treatises. All discussed the need to prove, the difficulty and the clarity of proofs. However, despite being a named mathematical object, proof was not yet constituted as an object of teaching. Geometry was *the* topic at stake and the ecological niche for proof to make sense.

It is with the emergence of state-based educational policies that in the mid-19th century the didactic transposition process began in earnest. However, at the very beginning of the century, Sylvestre-François Lacroix was already aware of the epistemic complexity of an educational project for mathematics. His Elements can be considered the first textbook as such: "it is his historical merit to have substantially contributed to the restructuring of a poorly-organized and scattered corpus of mathematical knowledge, guided by educational objectives" (Schubring, 1987, p. 43).

The didactic transposition process of Geometry continued to develop in order to respond to the growing needs of both industry and economy, and natural sciences as well. In this dynamic, mathematical proof remained untouched in its Euclidean norms. Organized as a professional body at the turn of the 20$^{th}$ century[46], mathematicians did recognize their responsibility towards mathematics education and contributed to the thinking on what mathematics should be as a content to be taught. Their concern was first to ensure that it kept its theoretical nature, not being reduced to a tool to the service of applications. They were aware of the problem of teaching students but they were not prepared to sacrifice their discipline; one of the first issues they considered was that of rigour. In fact, it was a question of deciding the acceptable limits of the didactic transposition of mathematical proof. This didactical process took another dimension when decision-makers published specification of the "knowledge to be taught". This took different forms in different nations depending on their educational organization and policy, but the movement was general.

---

[46] The community of professional mathematicians emerged as such with journals and societies in the course of the XIX° century. A first milestone is the creation of the French journal "Annales de Mathématiques Pures et Appliquées", founded and edited by J. D. Gergonne. It was published from 1810 to 1831. The first professional mathematical society is the Wiskundig Genootschap, founded in Amsterdam in 1778, but most others were founded in the second half of the nineteenth century (Bartle, 1995, p. 3).

Three key stages marked the spread of teaching of mathematics through the modern educational systems of the 20th century and the first decades of the 21st century. They determined the transposition of proof:

- *First stage*: The extension of the teaching of Geometry from grades 10-12 to grades 6-9; the latter being part of the compulsory school hence the teaching for all children. This required to choose a time for the introduction of mathematical proof. In general, the choice was to do it at grade 8. This was a real challenge and mostly a failure revealing the impossibility to escape *a rupture in the nature of the mathematics* "before" and an "after" the will to teach proof was made explicit.
- *Second stage*: The New math *coup de force* made mathematical proof a standard of validation in a unified mathematics. The rapid failure of this radical movement had the effect of the replacement in most of the following curricula of mathematical proof by deductive reasoning associated to a priority given to problem solving.
- *Third stage*: the willingness to introduce proof in the teaching at all grades of the compulsory school. This objective could not be reached without renouncing the standards of mathematical proof for teaching at the earlier grades. This comes with a vocabulary now including the words argumentation, justification and proof but without establishing a clear relation between the terms (and understanding the consequences of this absence on teaching).

The didactic transposition is a never-ending process as it is tightly dependant on the evolution of the society, its priority and shared epistemology, as well as on the evolution of mathematics itself and of the progress of knowledge on its teaching and learning.

## How to answer the question of truth before the availability of mathematical proof?

Until the beginning of the 20th century the idea of preparing the transition to mathematical proof preceding its introduction in the curricula was simply not considered. As it were, the difficulty of learning mathematical proof was recognized but viewed as the cost to pay for engaging in the learning of mathematics as a discipline. The failure of too many students in the context of the democratisation of mathematics education called for a response more effective than the one given when the psychology of development led to think that the transition could by itself be made possible by the children's access to the formal operational stage[47]. It calls for an evolution of the school epistemology necessary to answer the question that best captures the contemporary situation: *What should be taught before teaching mathematical proof?* Or, better: *how to answer the question of truth in the mathematics classroom before having mathematical proof available?*

The evolution of the didactic transposition is evidence of a pragmatic response from international bodies (e.g. TIMSS) and national educational institutions: teaching must allow the development of an argumentative competence – i.e. reasoned justification – as part of the early learning of mathematics, before the learning of mathematical proof[48]. For its part, research has undertaken work and projects to answer these questions [49]. But there is still a lack of contributions and results that are robust enough to allow curricula and teaching approaches to be designed in a reliable and efficient way. Contemporary research on proof and proving paints a complex picture of the relationship between argumentation and proof. There is debate and perhaps still a lack of consensus, although some points seem to be accepted:

- the structure of the text of an argumentation and of the text of an elementary proof are not radically foreign the one to the other.
- argumentation and proof have close relationships in the problem-solving process, but the transition from argumentation to proof needs a specific work.
- to get the status of a proof, an argumentation has to go through a social process ensuring its collective acceptance by the students *and* by the teacher

---

[47] The work of Jean Piaget had a historical and a significant impact on curriculum specification, possibly not only for what Constructivism brought about children learning, but because of the clarity of the Piagetian stages and their apparent simplicity. In practice, stakeholders and decision makers reduced the levels of argumentation to two categories: before and after the formal operational stage. The natural cognitive development was considered the determining factor of the levels of argumentation. The thought of Piaget was somewhat more sophisticated than that (see e.g. Piaget, 1973).

[48] i.e. the Euclidean standard of proof, which is still the reference structure of the mathematical discourse in classrooms.

[49] e.g. (Bieda et al., 2022; Hanna & de Villiers, 2012, Chapitre 15)

This closeness suggests that the question of the relation between argumentation and proof can have a response opening on a didactic transposition of proof in the form of an argumentation in the mathematics classroom acceptable from both a mathematical and a teaching perspective. However, argumentation has no formal status for the professional mathematician although it has a presence in the history of mathematics and in problem-solving processes. Then, with this absence of scientific status in contemporary mathematics: *could argumentation in the mathematics classroom be a didactic transposition of mathematical proof adjusted to the exigencies and constraints of teaching and learning at the compulsory school grades?*

Contemporary institutional texts, whether international or national, suggest a positive response, but do not share its foundations. What one gets instead, and in the first place what teachers get, is the idea of a seamless transition from arguing to proving in the mathematics classroom. Moreover, proof is present in the text of curricula and in their comment not as an object (i.e. content domain) but as a competence (i.e. cognitive domain). It induces that it cannot be taught directly but be stimulated and developed in situations that have mathematical and social characteristics to justify and give access to socio-mathematical norms (P. G. Herbst, 2002b; Yackel & Cobb, 1996). Research since the early 1970s has worked out characteristics of situations which engage students in establishing actively the validity of a statement (e.g. Brousseau, 1997, Chapitres 1, section 6) and the way that such situations challenge teachers (e.g. Ball, 1993; P. G. Herbst, 2002a; Lampert, 1990).

Andreas Stylianides (2007) proposed a characterization from which we could start:

> "Proof is a mathematical argument, a connected sequence of assertions for or against a mathematical claim, with the following characteristics:
>
> 1. it uses statements accepted by the classroom community (set of accepted statements) that are true and available without further justification;
> 2. it employs forms of reasoning (modes of argumentation) that are valid and known to, or within the conceptual reach of, the classroom community; and
> 3. it is communicated with forms of expression (modes of argument representation) that are appropriate and known to, or within the conceptual reach of, the classroom community." (ibid. p.291)

This characterization is appropriate, but it applies to any scientific discipline. It is too general, leaving open the main question for mathematics teachers and educators: *what would be the specific characteristics to add to account for the case of mathematics?*

Let us start from a remark: Mathematics develops on mathematics. This remark expresses the inward-looking epistemology which coins its form of abstraction. This does not contradict a mathematical activity which in many ways resembles the scientific activity, but as Christian Houzel (1979) put it: in mathematics the « already theorised knowledge… plays the role of the experimental instance »[50]. This is the origin of the radical abstractness of mathematics and of the specific nature of proof in this discipline.

The set of accepted statements – criterion 1 of Stylianides proposal – is more than a set, nor is it a repertoire: it is a set organised as *a system* which constitutes the material and the milieu[51] for the mathematical work. It was the objective of the construction of such a system which ultimately drove the writing of Euclid's Elements at the same time that it introduced a rupture with the sensory world[52]. The organization of the set of statements is the consequence of the fact that any of its elements is related to a subset of the whole by the links a proof establishes.

In the context of the classroom, this structured set of statements is not a proper theory insofar as its evolution is agile, including new admitted elements when necessary, and the modes of argumentation may vary in their nature having stronger roots in the community consensus than in a formalized ground. For this reason, I suggest to refer to it as a structured *Knowledge base*[53]. It will correspond to the first term *Theory* of the defining triplet of *Mathematical Theorem* in the sense of Alessandra Mariotti:

> "Proof is traditionally considered in itself, as if it were possible to isolate a proof from the statement to which it provides support, and from the theoretical frame within which this support makes sense. When one speaks of proof, all these elements, although not always mentioned, are actually involved at the same time, and it is not possible to grasp the sense of a mathematical proof without linking it to the other two elements: a statement and overall a theory." (2006, p. 183)

Moreover, we have to add two more constraints in order for an argumentation to reach the mathematical standard:

- on the one hand, that *a common norm of argumentation* is accepted and that any statement in the sequence of statements of an argumentation either is backed by an argumentation which meets the same requirements or comes from the knowledge base or,
- on the other hand, that it is ensured that any gap in the argumentation can be filled with an argumentation conforming the agreed norm.

This means establishing a practice that requires a deliberate transition from a *pragmatic* conception to a *rigorous* conception of proving. That is to say, a student shift from the position of a *practitioner* to the position of a *theoretician* (Balacheff, 1990).

Eventually, it is very unlikely that we will be able to find a solution for a seamless transition from arguing, in the general sense, to proving mathematically. For this reason, my position is to accept the creation of a didactical object: *mathematical argumentation*, and to work on its definition so that it provides a ground to building instructional bridges by creating the conditions for a socio-mathematical norm to become a precursor of mathematical proof.

This object, mathematical argumentation, cannot be conceived as a transposition of the mathematical proof unless one considers that the "social" function of the latter, within the scientific community, is constitutive of it (Balacheff, preprint). This would be an epistemological as well as a theoretical error: although being the product of a human activity that is certified at the end of a social process, mathematical proof is independent of a particular person or group (Delarivière et al., 2017). This will not be the case for a mathematical argumentation in the classroom. The standardisation of proof in mathematics, in addition to the institutional character of its reference (mathematical knowledge), has required its depersonalisation, its decontextualization and its timelessness. Yet argumentation is intrinsically carried by an agent, individual or collective, and dependent on the circumstances of its production.

---

[50] One could find cognitive consequences of this statement in (Tall et al., 2012)

[51] Milieu is used in the sense of the Theory of didactical situation (Brousseau, 1997).

[52] This is not contradictory with the use of mental experiments in some of the Elements' proofs, and with the recognition that the physical world and the other sciences contribute to the development of mathematic by the importation of certain intuitions, or by raising problems questioning mathematical concepts and models.

[53] (Balacheff, preprint)

The characteristics of mathematical argumentation must not only distinguish it from other types of argumentation used in scientific or non-scientific activities, in order to guarantee the possibility of the transition to the norm of mathematical proof, it must also be operational when it comes to arbitrating the students' proposals and eventually institutionalizing them in order to organize and to capitalize them in the classroom. Mathematical argumentation requires an institutionalisation. The recognition of its mathematical character cannot be reduced to a judgement on its form alone. How, for example, can we arbitrate the case of the generic example that balances the general and the particular, whose equilibrium is found at the end of a contradictory debate seeking an agreement that is as little as possible tainted by compromise?

Finally, proof is both a foundation and an organiser of knowledge. In the course of learning, it contributes to reinforcing knowledge evolution and to providing tools for its organisation. In teaching, it legitimises new knowledge and constitutes a system: knowledge and proof linked together provide the knowledge base with a structure which can work as a precursor to the theoretical ground mathematics need. The institutionalisation function of proof situations places explicit validation under the arbitration of the teacher who is ultimately the guarantor of its mathematical character. This social dimension, in the sense that scientific functioning depends on a constructed and accepted organisation, is at the heart of the difficulty of teaching proof in mathematics.

## Acknowledgements

I am grateful to Gert Schubring for his careful reading and attention to historical accuracy. I would particularly like to thank Evelyne Barbin, Gila Hanna, Patricio Herbst, Janine Rogalski and Nathalie Sainclair for their comments and suggestions. Of course, the responsibility for the entire paper, and especially for any remaining errors or misinterpretations, rests with me alone.